\documentclass[11pt]{article}

\usepackage[T1]{fontenc}
\usepackage[utf8]{inputenc}
\usepackage{lmodern}
\usepackage{amsmath,amssymb,amsthm}
\usepackage[margin=1in]{geometry}
\usepackage[hidelinks]{hyperref}
\hypersetup{
  pdftitle={Control of Oil Filtration Using Surface-Active Reagents},
  pdfauthor={Svetlana S. Mukhina},
  pdfkeywords={contact structure, contact transformation, Legendre surface,
    wave equation, oil filtration, Barenblatt model}
}

\allowdisplaybreaks

\newtheorem{theorem}{Theorem}

\newcommand{\R}{\mathbb{R}}
\newcommand{\M}{\mathcal{M}}
\newcommand{\A}{\mathcal{A}}
\newcommand{\Lg}{\mathcal{L}}
\newcommand{\vark}{\varkappa}
\newcommand{\dd}{\,\mathrm{d}}

\title{Legendre manifolds generated by the Barenblatt equations}

\author{Svetlana S. Mukhina\\
\small V.~A. Trapeznikov Institute of Control Sciences of the Russian Academy of Sciences,\\
\small Moscow, Russia\\
\small \texttt{ssmukhina@edu.hse.ru}}
\date{}

\begin{document}

\maketitle

\begin{abstract}
The paper is devoted to the construction of exact general solutions of the equations of oil filtration in a porous medium under the action of an active reagent injected into the reservoir. The choice of the reagent solution determines the control process of chemical oil flooding in hard-to-recover fields. The problem is described by a system of two first-order nonlinear partial differential equations (the Barenblatt model). Conditions are found on the Buckley--Leverett function, the adsorption function, and the distribution function of the reagent in oil under which the original system, by means of a nonlinear change of variables on a five-dimensional manifold endowed with a contact structure, is reduced to a linear wave system with constant coefficients. Six families of such functions are indicated. For them, contact transformations and Legendre surfaces corresponding to general solutions in implicit form are constructed.
\end{abstract}


\section*{Introduction}

In the development of oil fields, and primarily of fields with hard-to-recover reserves, chemical flooding is widely used. An aqueous solution of surface-active impurities is injected into an injection well; the solution moves through the oil reservoir and displaces oil toward a production well. The addition of the reagent changes the phase dynamics: the viscosities of water and oil, the densities, and the relative phase permeabilities change. As a result, the oil becomes more mobile, and it can be recovered even in cases where the energy of the reservoir and of classical flooding is no longer sufficient \cite{babalian}. In other words, the composition and concentration of the injected solution are control actions on the filtration process, and the choice of the reagent becomes a control problem.

The mathematical model of one-dimensional unsteady frontal displacement of oil by a solution of an active impurity in the large-scale approximation was proposed by G.~I.~Barenblatt \cite{barenblatt}. It has the form
\begin{equation}\label{eq:system}
\begin{cases}
s_t+H_x=0,\\[2pt]
\bigl(cs+\varphi(1-s)+a\bigr)_t+\bigl(cH+(1-H)\varphi\bigr)_x=0,
\end{cases}
\end{equation}
where $s=s(t,x)$ is the volume concentration of water in the porous space (water saturation), $c=c(t,x)$ is the concentration of the active reagent in water, $H=H(s,c)$ is the Buckley--Leverett function, $a=a(c)$ is the concentration of the deposited impurity on the pores of the reservoir, $\varphi=\varphi(c)$ is the concentration of the impurity in oil, $t$ is time, and $x$ is the spatial coordinate (the direction of filtration). The functions $\varphi(c)$, $a(c)$, and $H(s,c)$ depend on the physicochemical properties of the phases, the porous medium, and the reagent itself \cite{babalian}. System \eqref{eq:system} is of hyperbolic type \cite{derkach}.

The problem is to indicate reagent characteristics for which system \eqref{eq:system} admits exact solutions, and also to construct general solutions by the contactization method \cite{kushner2025}. Such solutions make it possible to estimate the influence of the reagent on the displacement front dynamics and can be used in choosing the injection regime. Results have been obtained for various forms of the functions $\varphi$, $a$, and $H$. Thus, in \cite{mukhina2023} the displacement by carbonated water was considered, when $\varphi(c)=\delta_1c+\delta_2$, $a(c)=0$, and the Buckley--Leverett function has the form $H=c(\mu s+\nu)$; a general multivalued solution was found. In \cite{mukhina2024}, the same model was studied in the framework of symplectic transformations. In contrast to \cite{mukhina2023,mukhina2024}, here the functions $\varphi$, $a$, and $H$ are considered in general form, and the admissible class of transformations is extended by the contactization method for first-order systems of equations. This extension is achieved by constructing geometric structures on a five-dimensional manifold, in particular, a contact form associated with the conservation law of the original system. It is precisely this extension of the class of transformations that makes it possible to describe a substantially wider set of reagents.

\section{The Barenblatt model and its geometrization}

With system \eqref{eq:system} in the space $\R^4$ with coordinates $t$, $x$, $s$, $c$, there is associated a pair of differential 2-forms \cite{kushner2007}:
\begin{equation}\label{eq:forms2}
\begin{aligned}
\omega_1&=H_s\,\dd t\wedge \dd s+H_c\,\dd t\wedge \dd c-\dd x\wedge \dd s,\\
\omega_2&=(c-\varphi)H_s\,\dd t\wedge \dd s
+\bigl((c-\varphi)H_c-(\varphi'-1)H+\varphi'\bigr)\dd t\wedge \dd c\\
&\quad-(c-\varphi)\,\dd x\wedge \dd s
+\bigl((s-1)\varphi'-a'-s\bigr)\dd x\wedge \dd c.
\end{aligned}
\end{equation}
The restrictions of these forms to the surface
$L=\{s=s(t,x),\,c=c(t,x)\}\subset\R^4$
vanish if and only if the pair of functions $(s(t,x),c(t,x))$ is a solution of system \eqref{eq:system}. System \eqref{eq:system} is a Jacobi system $E(\omega_1,\omega_2)$ \cite{kushner2007}, and solutions of system \eqref{eq:system} are identified with two-dimensional integral manifolds of the pair of forms \eqref{eq:forms2}.

Instead of the 2-forms \eqref{eq:forms2}, consider their linear combination
\[
\theta_1=a_{11}\omega_1+a_{12}\omega_2,\qquad
\theta_2=a_{21}\omega_1+a_{22}\omega_2
\]
with a nonzero determinant of the matrix $\|a_{ij}\|$. We choose the forms $\theta_1,\theta_2$ so that the conditions
\[
\theta_1\wedge\theta_2=0,\qquad
\theta_1\wedge\theta_1+\theta_2\wedge\theta_2=0
\]
hold. As a result, we obtain new 2-forms
\begin{equation}\label{eq:forms3}
\begin{aligned}
\theta_1&=-2H_c\sigma\zeta^{-1}\,\dd t\wedge \dd c
+2x\bigl(-\zeta H_{sc}+H_{ss}H_c\sigma\bigr)H_c^{-1}\zeta^{-1}\,\dd s\wedge \dd c\\
&\quad+\dd x\wedge \dd c-\dd t\wedge \dd s,\\
\theta_2&=\dd t\wedge \dd c-\dd x\wedge \dd s,
\end{aligned}
\end{equation}
where
\[
\sigma=(s-1)\varphi'-s-a',\qquad
\chi=(1-\varphi')H+\varphi',\qquad
\zeta=\sigma H_s+\chi.
\]
The form $\theta_2$ is closed. Consider the differential 1-form
$\theta=s\,\dd t+c\,\dd x$. Its exterior differential
$\dd\theta=\dd s\wedge \dd t+\dd c\wedge \dd x$ differs from the form
$\theta_2$ only in sign and vanishes on every solution of system
\eqref{eq:system}; that is, $\theta$ is a conservation law
\cite{kushner2007}.

Let us define the five-dimensional manifold
$\M=\R\times\R^4$ with coordinates $t$, $x$, $u$, $s$, $c$ and the projection
\[
\pi\colon\M\to\R^4,\qquad
\pi(t,x,u,s,c)=(t,x,s,c).
\]
The contact form on $\M$ is given by
\begin{equation}\label{eq:contact}
\vark=\dd u-\pi^*(\theta)=\dd u-s\,\dd t-c\,\dd x.
\end{equation}
The corresponding contact distribution $C=\ker\vark$ is generated by the vector fields
\begin{equation}\label{eq:fields}
\frac{\partial}{\partial t}+s\frac{\partial}{\partial u},\qquad
\frac{\partial}{\partial x}+c\frac{\partial}{\partial u},\qquad
\frac{\partial}{\partial s},\qquad
\frac{\partial}{\partial c}.
\end{equation}
The pair $(\pi^*(\theta_1),\vark)$ is called the contactization of the Jacobi system $E(\omega_1,\omega_2)$ \cite{kushner2025}. Note that the contact form $\vark$ is determined by the conservation law $\theta$ and depends on the system itself, in contrast to the Cartan form used in the geometry of differential equations.

\section{The contactization method}

A diffeomorphism $\Phi\colon\M\to\M$ is called contact if it preserves the contact structure:
$\Phi^*(\vark)=g\vark$ for some smooth function $g$ or, equivalently,
$\Phi^*(\vark)\wedge\vark=0$.

A submanifold $L\subset\M$ of dimension $2$ is called Legendre if
$\vark|_L=0$. If the Legendre submanifold $L$ is integral for the effective part of the form $\pi^*(\theta_1)$, then its projection $\pi(L)$ is a solution of the original system \cite{kushner2025}.

Recall that a 2-form $\gamma\in\Omega^2(\M)$ is called effective if it does not contain the differential $\dd u$ and is orthogonal to $\dd\vark$, i.e.
\[
\frac{\partial}{\partial u}\lrcorner\,\gamma=0,\qquad
\gamma\wedge \dd\vark=0.
\]
Let $D(C)$ be the module of vector fields lying in the contact distribution. Define the linear operator $\A\colon D(C)\to D(C)$ by
\[
X\lrcorner\,\pi^*(\theta_1)=\A X\lrcorner\,\dd\vark,\qquad X\in D(C).
\]
For the 2-forms from \eqref{eq:forms3}, it has the form
\[
\A=
\begin{pmatrix}
-1&0&0&\mu_2\\
\mu_1&1&\mu_2&0\\
0&0&-1&\mu_1\\
0&0&0&1
\end{pmatrix},
\]
where
\[
\mu_1=-2H_c\sigma\zeta^{-1},\qquad
\mu_2=2x\bigl(\sigma H_{sc}H_s-H_c\sigma H_{ss}+\chi H_{sc}\bigr)H_c^{-1}\zeta^{-1}
\]
in the basis \eqref{eq:fields}.

For a hyperbolic system $\A^2=1$; therefore, the eigenvalues of the operator are $\pm1$, and its eigenspaces form two two-dimensional characteristic distributions $C_+$ and $C_-$. Their first derivatives
\[
C_\pm^{(1)}=\langle X_\pm,Y_\pm,[X_\pm,Y_\pm]\rangle
\]
are three-dimensional, and their intersection in the five-dimensional manifold $\M$ occurs along the line
$\ell=C_+^{(1)}\cap C_-^{(1)}$. At each point $a\in\M$, the decomposition
\[
T_a\M=C_+(a)\oplus\ell(a)\oplus C_-(a)
\]
holds.

Let $P_1$, $P_2$, $P_3$ be the projections onto the distributions $D(C_+)$, $D(\ell)$, and $D(C_-)$, respectively. They have the form
\[
P_1=
\begin{pmatrix}
0&0&0&0&x\nu_1\\
\nu_2&\nu_3&-\nu_4&-x\nu_1&0\\
c\nu_2&c\nu_3&-c\nu_4&-cx\nu_1&sx\nu_1\\
0&0&0&0&-\nu_5\\
0&0&0&0&1
\end{pmatrix},
\]
\[
P_2=
\begin{pmatrix}
s\nu_1&c\nu_1&-\nu_1&0&0\\
s\nu_6&c\nu_6&-\nu_6&0&0\\
s\nu_6&c\nu_6&-\nu_6&0&0\\
0&0&0&0&0\\
0&0&0&0&0
\end{pmatrix},
\]
\[
P_3=
\begin{pmatrix}
\nu_7&-c\nu_1&\nu_1&0&-x\nu_1\\
-\nu_7(\sigma H_s+\chi)^{-1}&-c\nu_1H_c(\sigma H_s+\chi)^{-1}&-\nu_1H_c(\sigma H_s+\chi)^{-1}&x\nu_1&0\\
-\nu_7(\sigma H_s+\chi)^{-1}&-c\nu_1(\sigma H_s+\chi)^{-1}&\nu_1(\sigma H_s+\chi)^{-1}&xc\nu_1&-xs\nu_1\\
0&0&0&1&\nu_5\\
0&0&0&0&0
\end{pmatrix},
\]
where
\[
\nu_1=\frac{H_{sc}\zeta-H_{ss}\sigma H_c}{H_c\zeta},\qquad
\nu_2=\frac{(\varphi'+a')H_c}{\zeta},
\]
\[
\nu_3=\frac{\sigma H_s+\chi+c(\varphi'-1)H_c}{\zeta},\qquad
\nu_4=\frac{H_c(\varphi'-1)}{\zeta},
\]
\[
\nu_5=\frac{H_c\sigma}{\zeta},\qquad
\nu_6=\frac{(-\varphi'+1)\zeta-\sigma^2H_{ss}}{\zeta^2},\qquad
\nu_7=\frac{sH_{sc}(-\sigma H_s+\chi)+\sigma H_sH_c+H_c(\chi+s\sigma H_{ss})}{H_c\zeta}.
\]
Let us define the tensors \cite{kushner2025,kushner2022}
\[
q_{j,k}^{s}(X,Y)=-P_s\bigl([P_jX,P_kY]\bigr),\qquad
j,k,s=1,2,3,\quad s\ne j,k.
\]
The skew contraction operator of decomposable tensors $\alpha\otimes X$ and $\beta\otimes Y$ is defined by
\[
\langle\alpha\otimes X,\beta\otimes Y\rangle
=(Y\lrcorner\,\alpha)\wedge(X\lrcorner\,\beta).
\]
On indecomposable tensors it is defined by linearity. We define the 2-forms by the equalities \cite{kushner2008}
\[
\lambda_+=\langle q_{1,1}^{2},q_{2,3}^{1}\rangle,\qquad
\lambda_-=\langle q_{3,3}^{2},q_{1,2}^{3}\rangle.
\]

\begin{theorem}[\cite{kushner2025}]
A hyperbolic Jacobi system $E(\theta_1,\theta_2)$ for which the form $\dd\theta$ is nondegenerate is contact equivalent to a wave system if and only if
\[
\lambda_+=\lambda_-=0.
\]
\end{theorem}

If this condition is satisfied, then by the method of contact transformations the original Jacobi system with a nondegenerate conservation law can be reduced to a wave system with constant coefficients, and a Legendre surface corresponding to a multivalued solution can be constructed. For system \eqref{eq:system}, $\lambda_-=0$, and the form $\lambda_+$ has the form
\begin{equation*}\label{eq:lambda}
\begin{aligned}
\lambda_+={}&\Bigl(H_c\zeta^3H_{scc}-2H_c^2\zeta^2\sigma H_{ssc}
+H_c^3\zeta\sigma H_{sss}
+\bigl(3H_c^2\sigma^2H_{ss}-H_{cc}\zeta^2\bigr)\zeta H_{sc}\\
&-2H_c^2\Bigl(H_c\sigma^3H_{ss}
+\frac12\bigl(((1-H)a'+s-H)\varphi''-\sigma a''\bigr)\zeta\Bigr)H_{ss}\Bigr)
\,\dd s\wedge \dd c.
\end{aligned}
\end{equation*}
Thus, the condition for the reducibility of system \eqref{eq:system} to a wave system with constant coefficients is one nonlinear partial differential equation for the three functions $\varphi(c)$, $a(c)$, and $H(s,c)$. Let us indicate several particular cases in which $\lambda_+=0$:
\begin{enumerate}
\item $H(s,c)=(h_1s+h_2)r(c)$, the functions $a(c)$ and $\varphi(c)$ are arbitrary;
\item $H(s,c)=h_1s+h_2+r(c)$, the functions $a(c)$ and $\varphi(c)$ are arbitrary;
\item $H(s,c)=\bigl(\varphi'(s-1)-a'-s\bigr)r(c)+\dfrac{\varphi'}{\varphi'-1}$, the functions $a(c)$ and $\varphi(c)$ are arbitrary;
\item $H(s,c)=s\,r\!\left(\dfrac{c(s+1)}{s}\right)$, $a(c)=m_2$, $\varphi(c)=n_2$;
\item $H(s,c)=sr\!\left(c\exp\dfrac1s\right)$, $a(c)=m_2$, $\varphi(c)=n_2$;
\item $H(s,c)=s\,r(c)\pm\sqrt{(r(c)s-1)^2+4s}-1$, $a(c)=m_2$, $\varphi(c)=n_2$,
\end{enumerate}
where $h_1$, $h_2$, $m_2$, $n_2$ are arbitrary constants and $r(c)$ is an arbitrary function.

The first two conditions admit arbitrary functions $a$ and $\varphi$ and therefore describe the widest class of reagents. They generalize the cases considered earlier: for $\varphi(c)=\delta_1c+\delta_2$, $a(c)=0$, and $r(c)=c$, the first condition gives the Buckley--Leverett function $H=(h_1s+h_2)c$ from \cite{mukhina2023}, and the second condition with $H=h_1s+H_2(c)$ gives the function from \cite{mukhina2024}.

\section{Construction of exact solutions}

Let us turn to the construction of an exact solution of system \eqref{eq:system}. We choose the parameters $\varphi(c)$, $a(c)$, and $H(s,c)$ according to the condition $\lambda_\pm=0$.

\subsection{The first case}

Consider the Buckley--Leverett function of the form $H(s,c)=(h_1s+h_2)r(c)$. For simplicity, let the functions $a$ and $\varphi$ be linear:
\[
a(c)=m_1c+m_2,\qquad \varphi(c)=n_1c+n_2.
\]
Consider the five-dimensional manifold $\M$ with coordinates $t,x,u,s,c$ and the contact form \eqref{eq:contact}. Then the 2-form $\theta_1$ from \eqref{eq:forms3} has the form
\[
\theta_1=-\frac{2(h_1s+h_2)r'(c)(-n_1s+m_1+n_1+s)}
{(m_1h_1+n_1h_1+h_2n_1-h_2)r(c)-n_1}\,\dd t\wedge \dd c
+\frac{2h_1x}{h_1s+h_2}\,\dd s\wedge \dd c+\dd x\wedge \dd c-\dd t\wedge \dd s.
\]
Define the contact transformation $\Phi_1$:
\[
\Phi_1\colon\quad
\bigl\{t=h_1t,\quad x=q_2,\quad u=-h_2t+u,\quad
s=\frac{s-h_2}{h_1},\quad c=c\bigr\},
\]
then the 2-form $\omega=\Phi_1^*(\theta_1)$ has the form
\[
\omega=-\frac{2sr'(c)(-n_1s+h_1m_1+h_1n_1+h_2n_1-h_2+s)}
{(m_1h_1+n_1h_1+h_2n_1-h_2)r(c)-n_1}\,\dd t\wedge \dd c
-\frac{2x}{s}\,\dd s\wedge \dd c+\dd x\wedge \dd c-\dd t\wedge \dd s.
\]
The matrix representation of the operator $\A$ is
\[
\A=
\begin{pmatrix}
-1&0&0&\dfrac{2q_2}{p_1}\\[6pt]
-\dfrac{2p_1r'(p_2)K}{R}&1&-\dfrac{2q_2}{p_1}&0\\[6pt]
0&0&-1&-\dfrac{2p_1r'(p_2)K}{R}\\[6pt]
0&0&0&1
\end{pmatrix},
\]
where
\[
K=-n_1s+h_1m_1+h_1n_1+h_2n_1-h_2+s,\qquad
M=(m_1+n_1)h_1+h_2(n_1-1),\qquad
R=M\,r(p_2)-n_1.
\]
The eigenvalues of the operator are $\pm1$. The eigenvectors generate two two-dimensional contact characteristic distributions
\[
\begin{aligned}
V_1^+&=-\frac{Rx}{K\,r'(c)\,s^2}\frac{\partial}{\partial t}
-\frac{Rx}{sK\,r'(c)}\frac{\partial}{\partial u}
+\frac{\partial}{\partial s},\\
V_2^+&=\frac{R}{sK\,r'(c)}\frac{\partial}{\partial t}
+\frac{\partial}{\partial x}
+\left(c+\frac{R}{K\,r'(c)}\right)\frac{\partial}{\partial u},\\
V_1^-&=\frac{x}{s}\frac{\partial}{\partial t}
+x\frac{\partial}{\partial u}
-\frac{K\,r'(c)\,s}{R}\frac{\partial}{\partial s}
+\frac{\partial}{\partial c},\\
V_2^-&=\frac{\partial}{\partial x}+c\frac{\partial}{\partial u}.
\end{aligned}
\]
The intersection of the first derivatives $C_+^{(1)}\cap C_-^{(1)}$ is one-dimensional and is generated by the vector field
\[
Z=-\frac1s\frac{\partial}{\partial t}
-\frac{M\,r'(c)}{R}\left(\frac{\partial}{\partial x}+c\frac{\partial}{\partial u}\right).
\]
The first integrals of the distributions $C_-^{(1)}$ and $C_+^{(1)}$ have the form
\[
H_1^-=cx-u,\qquad H_2^-=\frac{K}{RM s},
\]
\[
H_1^+=c,\qquad
H_2^+=\frac{s\bigl(-t(n_1+m_1)h_1+(n_1-1)(-h_2t-cx+u)\bigr)r'(c)+xR}{sR}.
\]
Let us construct a contact transformation using the first integrals:
\[
\Phi_2\colon
\begin{cases}
t=cx-u,\\
x=c,\\
u=\dfrac{(-h_1t-h_2t-cx+u)n_1-h_1m_1t+h_2t+cx-u}{R},\\[6pt]
s=\dfrac{K}{sR},\\[6pt]
c=-\dfrac{M\bigl(s(-t(n_1+m_1)h_1+(n_1-1)(-h_2t-cx+u))r'(c)+xR\bigr)}{sR^2}.
\end{cases}
\]
Applying the transformation $\Phi_2^*$ to the form $\omega$ and extracting the effective part, we obtain
\[
\omega_{\mathrm{eff}}=\Phi_2^*(\omega)_{\mathrm{eff}}
=\frac{M\,r(c)-n_1}{s\bigl(M\,r(c)-n_1c+n_1-1\bigr)}
(\dd t\wedge \dd s-\dd x\wedge \dd c).
\]
The form $\omega_{\mathrm{eff}}$ corresponds to the form of the wave equation. The solution of the wave equation corresponds to the Legendre surface
\begin{equation}\label{eq:legendre}
\Lg=\{u=A(t)+B(x),\quad s=A'(t),\quad c=B'(x)\},
\end{equation}
where $A$ and $B$ are arbitrary functions. Applying the inverse transformations $\Phi_2,\Phi_1$ to the surface \eqref{eq:legendre}, we obtain the general solution of system \eqref{eq:system}.

Introduce the notation
\[
q_1=t,\qquad q_2=x,\qquad p_1=s,\qquad p_2=c,
\]
\[
C=(h_1+h_2)n_1+h_1m_1-h_2,\qquad D=Cr(p_2)-n_1,
\]
\[
Z=-p_2q_1(h_1p_1+h_2)r'(p_2)+h_2q_1r(p_2)+h_1u-h_2q_2,
\]
\[
\begin{aligned}
S&=-h_1q_1(m_1+n_1)r(p_2)
-p_2q_1(n_1-1)(h_1p_1+h_2)r'(p_2)\\
&\quad+(u+q_2)n_1+m_1q_2-u,
\end{aligned}
\]
\[
\begin{aligned}
Q&=-p_2q_1(n_1-1)(h_1p_1+h_2)r'(p_2)
+h_2q_1(n_1-1)r(p_2)\\
&\quad+(u-q_1+q_2)n_1+m_1q_2-u.
\end{aligned}
\]

\begin{theorem}
The Legendre surface corresponding to system \eqref{eq:system} under the conditions
$H(s,c)=(h_1s+h_2)r(c)$, $a(c)=m_1c+m_2$, $\varphi(c)=n_1c+n_2$ has the form
\[
\begin{cases}
\dfrac{S}{D}-A(Z)-B(-p_2)=0,\\[6pt]
A'(Z)+\dfrac{h_1\bigl((1-p_1)n_1+m_1+p_1\bigr)}{(h_1p_1+h_2)D}=0,\\[6pt]
B'(-p_2)-\dfrac{CQr'(p_2)}{D^2}=0.
\end{cases}
\]
Here $A$ and $B$ are arbitrary differentiable functions.
\end{theorem}

\subsection{The second case}

Consider the Buckley--Leverett function of the form $H(s,c)=h_1s+h_2+r(c)$. Then the 2-form $\theta_1$ from \eqref{eq:forms3} has the form
\[
\begin{aligned}
\theta_1={}&\dd x\wedge \dd c-\dd t\wedge \dd s\\
&-\frac{2r'(c)\bigl((s-1)\varphi'(c)-s-\alpha'(c)\bigr)\,\dd t}
{\bigl((s-1)\varphi'(c)-s-\alpha'(c)\bigr)h_1
+\bigl(1-\varphi'(c)\bigr)\bigl(h_1s+h_2+r(c)\bigr)+\varphi'(p_2)}
\wedge \dd c.
\end{aligned}
\]
In this subsection the adsorption function is denoted by $\alpha(c)$.

Let us define the five-dimensional manifold $\M$ with coordinates
$q_1=t$, $q_2=x$, $u=u$, $p_1=s$, $p_2=c$ and the given contact form \eqref{eq:contact}. The operator $\A$ has the form
\[
\A=
\begin{pmatrix}
-1&0&0&0\\
a_{21}&1&0&0\\
0&0&-1&a_{21}\\
0&0&0&1
\end{pmatrix},
\]
where
\[
a_{21}=
\frac{2r'(p_2)\bigl(\varphi'(p_2)p_1-\varphi'(p_2)-\alpha'(p_2)-p_1\bigr)}
{\varphi'(p_2)r(p_2)+\varphi'(p_2)h_1+\varphi'(p_2)h_2+\alpha'(p_2)h_1-\varphi'(p_2)-r(p_2)-h_2}.
\]
The eigenvectors generate two two-dimensional contact characteristic distributions
\[
V_1^+=
\frac{r'(p_2)\bigl((p_1-1)\varphi'(p_2)-p_1-\alpha'(p_2)\bigr)}
{\bigl((h_1+h_2+r(p_2)-1)\varphi'(p_2)+\alpha'(p_2)h_1-h_2-r(p_2)\bigr)}
\partial_{p_1}+\partial_{p_2},
\qquad
V_2^+=\partial_{q_2}+p_2\partial_u,
\]
\[
\begin{aligned}
V_1^-={}&-\frac{\bigl((h_1+h_2+r(p_2)-1)\varphi'(p_2)+\alpha'(p_2)h_1-h_2-r(p_2)\bigr)\partial_{q_1}}
{r'(p_2)\bigl((p_1-1)\varphi'(p_2)-p_1-\alpha'(p_2)\bigr)}
+\partial_{q_2}\\
&-\left(-p_2+
\frac{p_1\bigl((h_1+h_2+r(p_2)-1)\varphi'(p_2)+\alpha'(p_2)h_1-h_2-r(p_2)\bigr)}
{r'(p_2)\bigl((p_1-1)\varphi'(p_2)-p_1-\alpha'(p_2)\bigr)}
\right)\partial_u,
\end{aligned}
\]
\[
V_2^-=\partial_{p_1}.
\]
The intersection of the first derivatives of the contact characteristic distributions $C_+^{(1)}\cap C_-^{(1)}$ is generated by the vector field
\[
\begin{aligned}
Z={}&\frac{1-\varphi'(p_2)}
{(p_1-1)\varphi'(p_2)-p_1-\alpha'(p_2)}\partial_{q_2}\\
&+\frac{-p_2(\varphi'(p_2)-1)r'(p_2)+\alpha'(p_2)h_1
+\bigl(h_1+h_2+r(p_2)-1\bigr)\varphi'(p_2)-h_2-r(p_2)}
{r'(p_2)\bigl((p_1-1)\varphi'(p_2)-p_1-\alpha'(p_2)\bigr)}\partial_u.
\end{aligned}
\]
To simplify the finding of the first integrals, we choose the functions
\[
\alpha(p_2)=\alpha_1p_2^2,\qquad \varphi(p_2)=0.
\]
Let us construct the contact transformation $\Phi_1$:
\[
q_1^{\mathrm{new}}=q_1,
\]
\[
q_2^{\mathrm{new}}=
\frac{2Fq_1\alpha_1}{-\Delta}
\Bigl[r'J+\bigl((-p_2q_2+u)r'+q_2(2\alpha_1p_2h_1-h_2-r)\bigr)F
+2q_1r'\alpha_1p_2\Bigr],
\]
\[
\begin{aligned}
u^{\mathrm{new}}={}&\frac{1}{-\Delta}
\Bigl[2Fq_1\alpha_1(2\alpha_1p_2h_1+p_2r'-r-h_2)J\\
&+\bigl(p_2(-p_2q_2+u)r'+(2\alpha_1p_2h_1-h_2-r)u\bigr)F
+2r'\alpha_1p_2^2q_1\Bigr],
\end{aligned}
\]
\[
p_1^{\mathrm{new}}=F(2\alpha_1J+p_1),\qquad
p_2^{\mathrm{new}}=p_2,
\]
where
\[
\Delta=-2\alpha_1p_2h_1+r(p_2)+h_2,\qquad
F=\exp\left(\int\frac{r'(p_2)}{\Delta}\,\dd p_2\right),
\]
\[
J=\int\frac{p_2r'(p_2)}{\Delta}
\exp\left(-\int\frac{r'(p_2)}{\Delta}\,\dd p_2\right)\dd p_2.
\]
Applying the transformation $\Phi_1^*$ to the form $\theta_1$ and extracting the effective part, we obtain
\[
\omega_{\mathrm{eff}}=\Phi_1^*(\theta_1)_{\mathrm{eff}}
=-e^{\int\frac{r'(p_2)}{-2\alpha_1p_2h_1+r(p_2)+h_2}\,\dd p_2}
(\dd q_1\wedge \dd p_1-\dd q_2\wedge \dd p_2).
\]
This form $\omega_{\mathrm{eff}}$ corresponds to the form of the wave equation. Let us introduce the notation
\[
q_1=t,\quad q_2=x,\quad p_1=s,\quad p_2=c,
\]
\[
F(p_2)=-2\alpha_1p_2h_1+r(p_2)+h_2,\qquad
G(p_2)=2\alpha_1p_2h_1-h_2-r(p_2),
\]
\[
E(p_2)=\exp\left(-\int\frac{r'(p_2)}{F(p_2)}\,\dd p_2\right),\qquad
J(p_2)=\int\frac{p_2r'(p_2)E(p_2)}{F(p_2)}\,\dd p_2,\qquad
z=-h_1q_1+q_2,
\]
\[
W=\frac{r'(p_2)\bigl(2\alpha_1p_2(h_1q_1-q_2)-p_2q_1r'(p_2)+u-q_1G(p_2)\bigr)}{G(p_2)},
\]
\[
R=-p_2\bigl(2\alpha_1p_2(h_1q_1-q_2)-p_2q_1r'(p_2)+u\bigr)r'(p_2)+uG(p_2).
\]

\begin{theorem}
The Legendre surface corresponding to system \eqref{eq:system} under the conditions
$H(s,c)=h_1s+h_2+r(c)$, $\alpha(p_2)=\alpha_1p_2^2$, $\varphi(p_2)=0$ has the form
\[
\begin{cases}
\dfrac{ER}{G}-B(EW)-A(z)-2\alpha_1(h_1q_1-q_2)J=0,\\[6pt]
-2\alpha_1J+p_1E+A'(z)=0,\\[4pt]
p_2+B'(EW)=0.
\end{cases}
\]
Here $A$ and $B$ are arbitrary differentiable functions.
\end{theorem}

\subsection{The third case}

Consider the Buckley--Leverett function, the adsorption function, and the concentration of the impurity in oil of the form
\[
H(p_1,p_2)=p_1\,r\!\left(\frac{p_2(p_1+1)}{p_1}\right),\qquad
a(p_2)=m_2,\qquad \varphi(p_2)=n_2.
\]
Then the 2-form $\theta_1$ from \eqref{eq:forms3} has the form
\begin{equation}\label{eq:theta7}
\theta_1=-\dd q_1\wedge \dd p_1
+\frac{2(p_1+1)p_1}{p_2}\,\dd q_1\wedge \dd p_2
+\dd q_2\wedge \dd p_2
-\frac{2q_2}{p_1+1}\,\dd p_1\wedge \dd p_2.
\end{equation}
Apply the contact transformation
\[
\Phi=\{q_1=q_1,\ q_2=q_2,\ u=u-q_1,\ p_1=p_1-1,\ p_2=p_2\}
\]
to the 2-form $\theta_1$ from \eqref{eq:theta7}:
\[
\omega=\Phi^*(\theta_1)
=-\dd q_1\wedge \dd p_1
+\frac{2(p_1-1)p_1}{p_2}\,\dd q_1\wedge \dd p_2
+\dd q_2\wedge \dd p_2
-\frac{2q_2}{p_1}\,\dd p_1\wedge \dd p_2.
\]
The operator $\A$ has the form
\[
\A=
\begin{pmatrix}
-1&0&0&\dfrac{2q_2}{p_1}\\[6pt]
\dfrac{2p_1(p_1-1)}{p_2}&1&-\dfrac{2q_2}{p_1}&0\\[6pt]
0&0&-1&\dfrac{2p_1(p_1-1)}{p_2}\\[6pt]
0&0&0&1
\end{pmatrix}.
\]
The contact characteristic distributions have the form
\[
V_1^+=\frac{q_2p_2}{p_1^2(p_1-1)}\partial_{q_1}
+\frac{q_2p_2}{p_1(p_1-1)}\partial_u+\partial_{p_1},
\]
\[
V_2^+=-\frac{p_2}{p_1(p_1-1)}\partial_{q_1}
+\partial_{q_2}
-\left(-p_2+\frac{p_2}{p_1-1}\right)\partial_u,
\]
\[
V_1^-=\frac{q_2}{p_1}\partial_{q_1}+q_2\partial_u
+\frac{p_1(p_1-1)}{p_2}\partial_{p_1}+\partial_{p_2},
\]
\[
V_2^-=\partial_{q_2}+p_2\partial_u.
\]
The intersection of the first derivatives of the characteristic distributions is generated by the vector field
\[
Z=-\frac{1}{p_1}\partial_{q_1}-\frac{1}{p_2}\partial_{q_2}-\partial_u.
\]
The first integrals of the derivatives of the characteristic distributions determine the contact transformation
\[
\begin{aligned}
\Phi_1=\Bigl\{&q_1=p_1^2q_1-up_1-q_2,\quad
q_2=\frac{2p_1q_1-u}{p_2p_1+1},\\
&u=\frac{2p_1^2q_1-p_1p_2q_2-up_1-q_2}{p_2p_1+1},\quad
p_1=\frac{1}{p_2p_1+1},\quad p_2=p_1\Bigr\}.
\end{aligned}
\]
Applying the transformation $\Phi_1^*$ to the 2-form $\omega$, we obtain
\[
\omega_{\mathrm{eff}}=\Phi_1^*(\omega)_{\mathrm{eff}}
=\frac{p_1}{p_1p_2+1}(\dd q_1\wedge \dd p_1-\dd q_2\wedge \dd p_2).
\]
This form $\omega_{\mathrm{eff}}$ corresponds to the form of the wave equation.

\begin{theorem}
The Legendre surface corresponding to system \eqref{eq:system} under the conditions
\[
H(p_1,p_2)=p_1r\!\left(\frac{p_2(p_1+1)}{p_1}\right),\qquad
a(p_2)=m_2,\qquad \varphi(p_2)=n_2
\]
has the form
\[
\begin{aligned}
&\frac{2(u+q_1)}{p_2}-\frac{2q_1}{p_2}+\frac{q_2}{p_1+1}-2q_2
-A\left(\frac{q_2p_1p_2-u(p_1+1)}{p_2^2(p_1+1)}\right)
-B(-p_2q_2+q_1+u)=0,\\
&p_2+A'\left(\frac{q_2p_1p_2-u(p_1+1)}{p_2^2(p_1+1)}\right)=0,\\
&-\frac{p_1}{(p_1+1)p_2}+B'(-p_2q_2+q_1+u)=0.
\end{aligned}
\]
Here $q_1=t$, $q_2=x$, $p_1=s$, $p_2=c$, and $A$ and $B$ are arbitrary differentiable functions.
\end{theorem}

\subsection{The fourth case}

Consider the fifth condition
\[
H(s,c)=sr\!\left(c\exp\frac1s\right),\qquad a(c)=m_2,\qquad \varphi(c)=n_2.
\]
Then the 2-form $\theta_1$ has the form
\begin{equation}\label{eq:theta8}
\theta_1=-\dd q_1\wedge \dd p_1
+\frac{2p_1^2}{p_2}\,\dd q_1\wedge \dd p_2
+\dd q_2\wedge \dd p_2
-\frac{2q_2}{p_1}\,\dd p_1\wedge \dd p_2.
\end{equation}
The operator $\A$ has the form
\[
\A=
\begin{pmatrix}
-1&0&0&\dfrac{2q_2}{p_1}\\[6pt]
\dfrac{2p_1^2}{p_2}&1&-\dfrac{2q_2}{p_1}&0\\[6pt]
0&0&-1&\dfrac{2p_1^2}{p_2}\\[6pt]
0&0&0&1
\end{pmatrix}.
\]
The eigenvectors generate two two-dimensional contact characteristic distributions
\[
V_1^+=\frac{q_2}{p_1}\partial_{q_1}+q_2\partial_u
+\frac{p_1^2}{p_2}\partial_{p_1}+\partial_{p_2},
\qquad
V_2^+=\partial_{q_2}+p_2\partial_u,
\]
\[
V_1^-=\frac{q_2p_2}{p_1^3}\partial_{q_1}
+\frac{q_2p_2}{p_1^2}\partial_u+\partial_{p_1},
\qquad
V_2^-=-\frac{p_2}{p_1^2}\partial_{q_1}+\partial_{q_2}
+\left(p_2-\frac{p_2}{p_1}\right)\partial_u.
\]
The intersection of the first derivatives of the contact characteristic distributions $C_+^{(1)}\cap C_-^{(1)}$ is generated by the vector field
\[
Z=\frac{p_2}{p_1^3}\partial_{q_1}.
\]
The first integrals of the derivatives of the characteristic distributions determine the contact transformation $\Phi$:
\[
\begin{aligned}
\Phi=\Bigl\{&q_1=q_2\ln p_1+p_1q_1-u,\quad
q_2=-\frac{p_1q_1+q_2}{p_1(\ln p_1-p_2)},\\
&u=-\frac{q_2\ln p_1+p_1q_1-p_2q_2+q_2}{\ln p_1-p_2},\quad
p_1=-\frac{1}{\ln p_1-p_2},\quad p_2=p_1\Bigr\}.
\end{aligned}
\]
Then the 2-form $\theta_1$ from \eqref{eq:theta8} has the form
\[
\omega_{\mathrm{eff}}=\Phi(\theta_1)_{\mathrm{eff}}
=\frac{1}{-\ln p_1+p_2}(\dd q_1\wedge \dd p_1-\dd q_2\wedge \dd p_2).
\]
This form $\omega_{\mathrm{eff}}$ corresponds to the form of the wave equation.

\begin{theorem}
The Legendre surface corresponding to system \eqref{eq:system} under the conditions
$H(s,c)=sr\!\left(c\exp\dfrac1s\right)$, $a(c)=m_2$, $\varphi(c)=n_2$ has the form
\[
\begin{aligned}
&-\frac{1}{p_1}\Biggl(-A\left(\frac{(p_2q_2-u)p_1-p_2q_2}{p_1p_2}\right)p_1
-B(-p_2q_2+u)p_1\\
&\qquad-p_1(-p_2q_2+u)\ln(p_2q_2)
+p_1(-p_2q_2+u)\ln(q_2)
+(-p_2q_2-q_1+u)p_1+p_2q_2\Biggr)=0,\\
&p_2+A'\left(\frac{(p_2q_2-u)p_1-p_2q_2}{p_1p_2}\right)=0,\\
&\ln(p_2q_2)+\frac{-p_1\ln(q_2)+1}{p_1}+B'(-p_2q_2+u)=0.
\end{aligned}
\]
Here $q_1=t$, $q_2=x$, $p_1=s$, $p_2=c$, and $A$ and $B$ are arbitrary differentiable functions.
\end{theorem}

\subsection{Verification of the solution}
Let us verify the solution found. Choose the functions from \eqref{eq:legendre}
$A(x)=x^2$, $B(x)=x^2$; then, projecting the Legendre surface from $\M$ onto $\R^4$, we obtain the solution of the initial system \eqref{eq:system} in parametric form
\[
S\colon
\begin{cases}
t=\dfrac{1+\tau\gamma^2-\tau\ln\left(\dfrac1\gamma\right)}
{2\exp\left(\dfrac1\tau\right)\gamma\tau\,r'\!\left(\gamma\exp\left(\dfrac1\tau\right)\right)},\\[12pt]
x=r\!\left(\gamma\exp\left(\dfrac1\tau\right)\right)t(\tau,\gamma)
-\dfrac{\left(\tau\ln\left(\dfrac1\gamma\right)-1\right)^2+\tau\gamma^2(\tau+2)}{4\tau^2},\\[12pt]
s=\tau,\qquad c=\gamma,
\end{cases}
\]
where $\gamma,\tau$ are parameters. Indeed,
$\omega_1|_S=\omega_2|_S=0$.

\section*{Conclusion}

The paper considered the problem of one-dimensional oil displacement by a solution of an active reagent. Necessary and sufficient conditions are obtained under which it is contact equivalent to a wave system with constant coefficients. The conditions are formulated as one nonlinear equation for the Buckley--Leverett function, the adsorption function, and the distribution function of the reagent in oil. Six families of its solutions are described. For four cases, contact transformations are constructed, characteristic distributions and their first integrals are found, and the form is reduced to the wave form; a Legendre surface corresponding to the general solution of the system is obtained. The constructed solution is verified. By choosing initial concentration distributions, one can solve the Cauchy problem and visualize the solutions.

\section*{Acknowledgments}

This work was supported by the Foundation for the Development of Theoretical Physics and Mathematics ``BASIS'', grant No.~23-7-5-52-1.

\end{document}